\documentclass[a4paper, final, 12pt]{article}

\usepackage{eurosym}
\usepackage{amsmath, amssymb, latexsym, amscd, amsthm,amsfonts,amstext}
\usepackage[mathscr]{eucal}
\usepackage{graphicx}
\usepackage{subfig}
\usepackage{float}
\usepackage{color}
\usepackage{hyperref}
\usepackage[utf8]{inputenc}
\usepackage[english]{babel}

\usepackage{siunitx}

\numberwithin{equation}{section}
\newcommand{\appropto}{\mathrel{\vcenter{
  \offinterlineskip\halign{\hfil$##$\cr
    \propto\cr\noalign{\kern1pt}\;\sim\cr\noalign{\kern-2pt}}}}}

\begin{document}

\title{ A Game of Simulation:\\
Modeling and Analyzing the Dragons \\
of Game of Thrones}
\author{Zheng Cao, Brody Bottrell, Jiayi Gao, Mark Pock, Vinsensius  
\and University of Washington, Seattle, USA
\and Department of Mathematics
\and zc68@uw.edu
\and Department of Astronomy
\and brodyb03@uw.edu
\and Academy for Young Scholars
\and jerrygao@uw.edu
\and Department of Computer Science
\and markpock@uw.edu 
\and Department of Applied Mathematics
\and vxvinsen@uw.edu
}
\date{}
\maketitle

\begin{abstract}

This paper outlines two approaches for mathematical modelling and analysis of hypothetical creatures -- in particular, the dragons of HBO's television series Game of Thrones (GOT). Our first approach, the forward model (FM), utilises quasi-empirical observations of various features of GOT dragons. We then mathematically derive the growth rate, other dimensions, energy consumption, etc. In the backward model (BM), we use projected energy consumption by given ecological impact to model an expected dragon in terms of physical features. We compare and contrast both models to examine the plausibility of a real-world existence for our titular dragons and provide brief analyses of potential impacts on ecology.

\end{abstract}

\vspace{0.5em}
\textbf{\text{Keywords:}}

Numerical Analysis, Growth Modeling, 3D Modeling, Energy Metabolism, Dragon

\newpage

\vspace{0.5em}
\tableofcontents

\newpage

\section{Introduction}
\hspace{1em}
The dragon occupies a high place in our collective mythology, cast as a primordial, immense serpent possessing vast power. If the dragon is benevolent, it uses its power to guide kings, establish nations, and raise up heroes. It comes to embody our precarious yet dependent relationship to nature. A malevolent dragon, though, takes on the role of an ur-evil, and one that stems from a deep place in the human psyche. Its power instead becomes an oppressive terror, a living symbol of eternal, tyrannical, unshakable authority.

It is no accident, then, that dragons occupy an equally large place in our high fantasy. The archetypal epic fantasy saga in Tolkien's Middle-Earth created the general picture of dragons in the popular imagination. Like Smaug, they are four-legged, two-winged beasts of a vaguely lizard-like (or sometimes more serpentine) form whose fiery breath burns hot enough to incinerate cities.

Probably the cornerstone name of modern fantasy is Game of Thrones, the hit HBO adaptation of George R. R. Martin's \textit{A Song of Ice and Fire}. Its characters -- Ned Stark, Jon Snow, Daenerys Targaryen -- have become household names. In spite of the widely negative critical reaction to the final season, the franchise only continues to grow, with several more prequel and sequel shows in production, most notably House of the Dragon, which seeks to chronicle the massive civil war occurring some two hundred years before the events of Game of Thrones.

We intend -- if the reader will bear the exercise and suspend their disbelief -- to mathematically model these dragons, these fearsome beasts of war, in terms of growth, dimensions, and energy consumption. As a notional guiding question, we seek to understand whether the dragon of the modern popular imagination (as represented by the dragons of Game of Thrones) has any place in reality. We have attempted to be as sound as possible both in mathematical argument and empirical observation (as much as either can be achieved under the circumstances).

\subsection{Core Questions}
\hspace{1em}
This paper will explore the plausibility of dragons' existence in the real world via mathematical argument. We then attempt to probe the sustainability of their existence. We thus begin by posing a series of core questions as follows.

\begin{enumerate}
    \item What are the most salient features of the dragons of Game of Thrones?
    \item How would the existence of dragons impact environments?
    \item How much must dragons consume to support their outlined behaviors?
    \item How much land area in a typical environment is required to support a dragon? 
\end{enumerate}

\section{Strategy}
\hspace{1em}
To verify or falsify the plausibility of the existence of dragons, we attempt two strategies, a `forward model' and a `backward model.'

\subsection{Forward Model}
\hspace{1em}
We begin by laying out a cohesive slate of useful dimensions for which data can be feasibly gathered. Our first priority becomes to find scenes where Drogon and Daenerys are roughly the same distance away from the camera. Having vetted said scenes, we make the assumption for the sake of data collection that measurements taken from the relevant shots are all (approximately) proportional to the actual figures. Under this assumption, we are able to use the actual height of the actors to extrapolate approximate measurements.

We then attempt to model the volume using 3D modeling techniques. We extrapolate the volume to potential energy consumption via an extension of metabolic rate calculations for birds. Finally, using these energy consumption figures, we provide an overview of the type of ecosystem and amount of land required to sustain a dragon and allow it to exist in the real world.

All of the measurements we collect from the TV series are with respect to Drogon, as Drogon is the only dragon that developed without artificially imposed constraints (unlike Drogon's two siblings, confined at a young age). Drogon is also the most prominent of the three, so we were most able to gather information.

\subsection{Backward Model}
\hspace{1em}
Starting from conservation of energy
based on the literal texts of how much food a dragon consumed a day
to the amount of energy a dragon needed to sustain its life (based on real-life creatures, such as Komodo dragon)
then to compute an ideal dragon for which all its physical habits' energy requirements can be met, for example, flying, hunting, mating, fire-breathing, etc.)

We then attempted to model, analyze, and distribute the amount of of energy required for stated behaviours -- e.g. flight and fire breath -- to test whether such a dragon can exist in the real world.

\section{Notation}
\begin{enumerate}
    \item [$t \;\;$] Number of elapsed years since birth of dragon
    \item [$L(t) \;\;$] Length of dragon at time $t$ in meters
    \item [$L_H(t) \;\;$] Length of dragon's head/snout at time $t$ in meters
    \item [$A \;\;$] Maximum (`adult') length constant in meters
    \item [$l \;\;$] Distance from a lateral cross-section perpendicular to the spine to the end of the snout, measured along the spine in meters
    \item [$W(l) \;\;$] Longest straight-line distance between the center of the cross-section at $l$ and a point on its boundary in meters
    \item [$W(t) \;\;$] `Average width' of dragon at time $t$ in meters
    \item [$C(l) \;\;$] Cross-sectional area of the dragon at $l$ in square meters
    \item [$H(l) \;\;$] Height of the dragon at $l$ in meters
    \item [$H(t) \;\;$] `Average height' of dragon at time $t$ in meters
    \item [$V_T(t) \;\;$] Volume of the torso of the dragon at time $t$ in cubic meters
    \item [$V(t) \;\;$] Volume of dragon at time $t$ in cubic meters
    \item [$\rho \;\;$] Average density constant in kilograms per cubic meter
    \item [$M(t) \;\;$] Mass of dragon at time $t$ in kilograms
    \item [$R_B(t) \;\;$] Basal metabolic rate of dragon at time $t$ in kilocalories per hour
    \item [$R_S(t) \;\;$] Standard metabolic rate of dragon at time $t$ in kilocalories per hour
    \item [$R_F(t) \;\;$] Metabolic rate during flight of dragon at time $t$ in kilocalories per hour
    \item [$P_F \;\;$] Percentage of time spent flying (dimensionless)
    \item [$C(t) \;\;$] Daily energetic consumption not accounting for fire breathing in kilocalories per day
\end{enumerate}

\section{Data Collection}
\hspace{1em}
For the data collection portion, we first obtained measurements of people like Jon Snow and Daenerys Targaryen, especially their head height, to serve as the known variable when scaling to obtain measurements. With this information we chose clips from the show that included both of the faces of the dragons and Daenerys/Jon were in shot. We used the website: \url{https://eleif.net/} to check the scaling between a dragon’s head height to Deanery/Jon. We did this multiple times per season in order to get a range to allow for error on the measurements. We repeated this process to get the head width as well. With the head height and width, we no longer needed Daenerys or Jon in a shot to get new measurements, so by then it was onto getting the dragon’s full-body length. How we approached this was by again, getting specific shots of the dragons, but this time with their full body out, so we could scale how big the head it compared to the rest of the body, with the same tool as before. We would use the uncertainty ranges from before in order to properly get some kind of uncertainty on the body length as well.

\begin{center}
\hspace{20pt}
\includegraphics[scale=0.3]{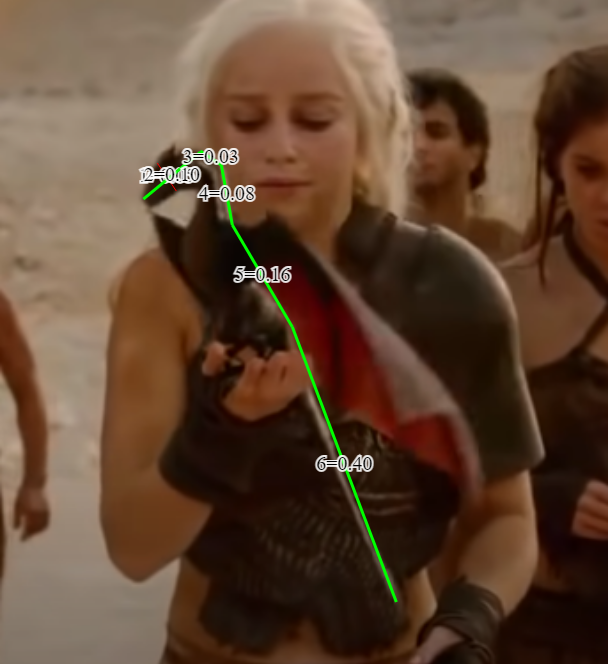}\\
Figure 1: Season 2 Example Data Image 1 \cite{GoT}
\end{center}

What we see in the table below is all the information from each season by using this measuring method.

\bigskip \vspace{0.75cm} 
\textbf{Table 1. Dragon Data 1.0}
\label{dragondatatab1}
\begin{center}
\begin{tabular}{|l|l|l|l|l|l|l|}
\hline
Age(yrs) & Head Height(m) & Head Length(m) & Body Width(m) & Body Length(m) \\ \hline
0 & 0.046 & 0.0575 & 0.050 & 0.436 \\ \hline
0.5 & 0.034-0.046 & 0.044-0.059 & 0.037-0.050 & 0.520... \\ \hline
2 & 0.122-0.133 & 0.138 & 0.132-0.144 & 1.65... \\ \hline
3 & 0.253 & 0.0575 & 0.275 & 4.408... \\ \hline
4 & 0.6 (0.57) & 0.735 & 0.619 & 12.988... \\ \hline
5 & 1.03-1.076 & 1.65-1.719 & 1.12-1.78 & 30.59... \\ \hline
6 & 2 (2.04) & 2.16-2.57 & 2.22 & 37.08... \\ \hline
7 & 2 (2.04) & 2.16-2.57 & 2.22 & 53.26... \\ \hline
\end{tabular}
\end{center}

What differs between the first and second table, is the certainty present. Table 2 provides us the best data collected from table 1, expertly showing the growth of every dimension on the dragon.

\bigskip \vspace{0.75cm} 
\textbf{Table 2. Dragon Data 2.0}
\label{dragondatatab2}
\begin{center}
\begin{tabular}{|l|l|l|l|l|l|l|}
\hline
Age(yrs) & Head Height(m) & Head Length(m) & Body Width(m) & Body Length(m) \\ \hline
0 & 0.046 & 0.0575 & 0.050 & 0.436... \\ \hline
0.5 & 0.046 & 0.059 & 0.050 & 0.956... \\ \hline
2 & .128 & 0.138 & .138 & 2.900... \\ \hline
3 & 0.253 & 0.275 & 0.275 & 6.396... \\ \hline
4 & 0.57 & 0.735 & 0.619 & 13.953... \\ \hline
5 & 1.053 & 1.685 & 1.45 & 31.239... \\ \hline
6 & 2.04 & 2.365 & 2.22 & 40.748... \\ \hline
7 & 2.04 & 2.57 & 2.22 & 49.758... \\ \hline
\end{tabular}
\end{center}

The next data collection that we will do is to do volume data collection to calculate the mass of the dragon. To do that, we will be using open-source 3D models of dragon available online. Unfortunately, we are not able to find Drogon's model that is free to public. So, we will be using other dragon model that is similar to Drogon in appearance. We are using an open-source software called Blender to calculate the volume of the model. To calculate the model, we will first scale the model such that it has the same body length as shown in table \ref{dragondatatab1}. Then, Blender will calculate the volume of the model so that it can be used for calculating the mass of the dragon. This process will be further explained in section 4.4: Building the Model.

\begin{center}
\hspace{20pt}
\includegraphics[scale=0.2]{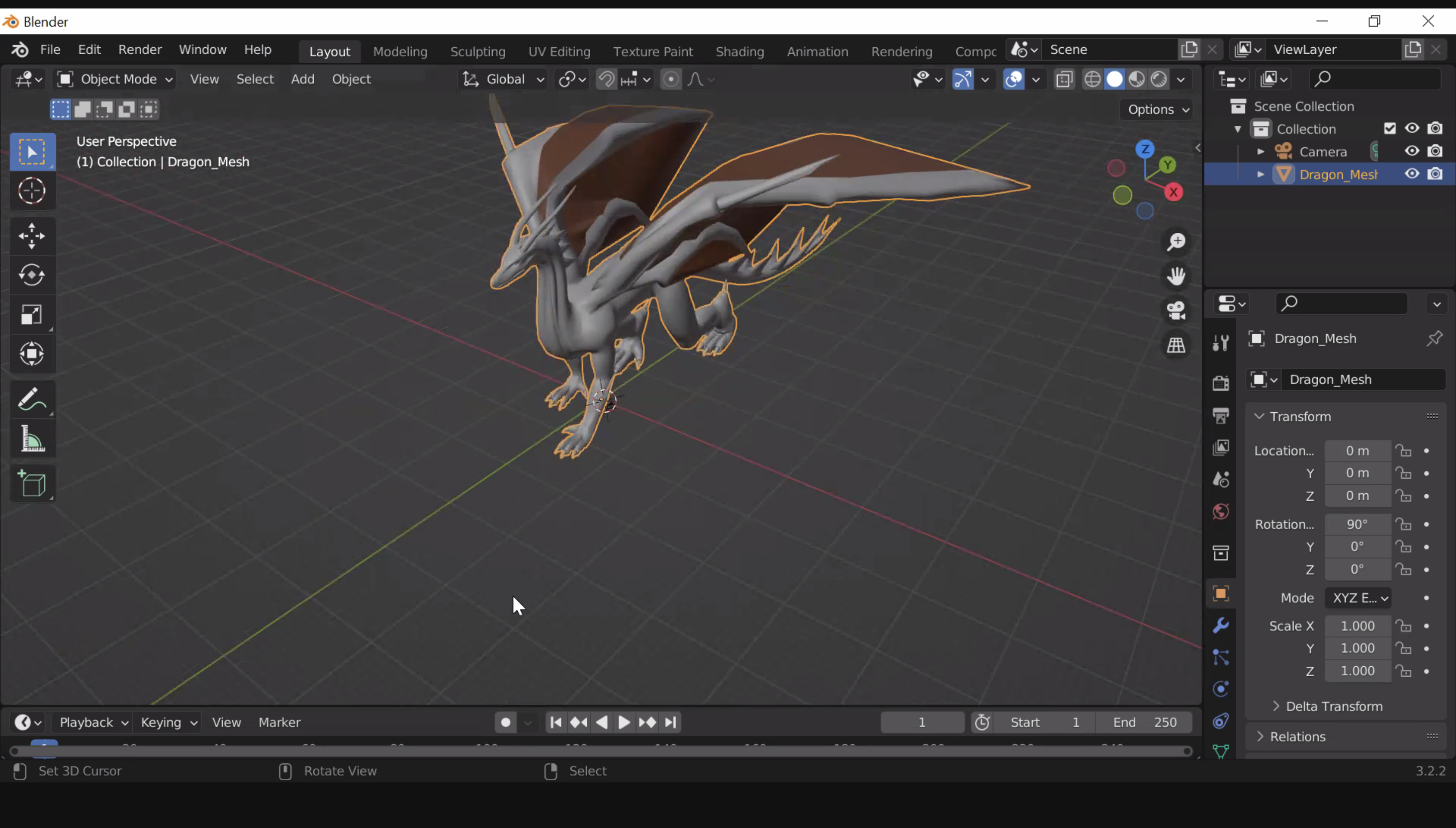}\\
Figure 2: 3D Modeling Example 1
\end{center}

\section{Forward Modeling}
\hspace{1em}
The section has three main objectives. First, we generalize the growth over time of a dragon through the use of mathematical modelling, drawing directly on scenes from \textit{Game of Thrones} and using the proscribed extrapolation techniques. Second, we extend this growth function $(L(t))$ to find an expression for the mass $M(t)$ under the assumptions given in 5.1. Finally, we extend the mass to a proposed energetic consumption and use it to extrapolate a daily intake, speculating on its ecological impacts.

\subsection{Assumptions}
\hspace{1em}
We frame Section 5.2 in terms of four assumptions from which can be derived an approximate expression for the mass of the dragon in terms of its length.

\begin{enumerate}

    \item We assume that height is proportional to width. While this is not entirely grounded, it is rather necessary to proceed -- we found it difficult to obtain accurate measurements of long-term change in height.
    \[H(t) \propto W(t) \]
    
    \item We assume that the average cross-sectional area of the dragon grows proportionally to the product of its width and height over time. Essentially, we assume that the \emph{average} shape of the cross-sectional area is roughly the same no matter the age of the dragon, or that any changes in one part of the body are roughly cancelled out by changes in another part of the body. Since it is unfeasible to measure the exact shape of the cross-sectional area for all parts of the body over all given timeframes, we have to make do with the assumption as follows.
    \[\dfrac{1}{L(t)}\int_{0}^{L(t)}C(l)dl \propto W(t) \cdot H(t)\]
    
    \item Most flying creatures have very light body weights with respect to their often much larger volumes -- a requirement for floating and flying. Therefore, we have assumed that the average density of a dragon will be less than or equal to the density of pure water, or $997$ kg/m$^2$.
    
    \[\rho \leq 997 \text{ kg/m}^2\]
    
    \item We take the wings of the dragon to  account for a negligible amount of mass when compared against the remainder of the body, assuming that most of the dragon's mass is in its torso. This does not necessarily entail that the torso of the dragon comprises the majority of its volume but rather should be taken to mean that the density is highly variable and that our chosen $\rho$ value acts as a good approximation by which to incorporate the limited mass of the low-density wings.
    
    \[\rho\int_{0}^{L(t)}C(l)dl \approx M(t)\]
\end{enumerate}

\subsection{Dimensions}
\hspace{1em}
Based on the data we collected from the literature and from pixel measurement, we conclude, that as is generally expected for animal growth, dragons mostly exhibit long-term asymptotic growth up to a distribution around final adult age. We do however, note that the oldest surviving specimen, Balerion, was said not to have `stopped growing' until two or three years before his death at around three hundred years of age -- an ambiguous description in some sense, since we are unable to pin to which dimension it refers. The assumption is thus that the descriptor refers to either wingspan or width, and that if any growth in length occurred it was minimal at best.

As such, our priority becomes finding a suitable model for long-term asymptotic growth. While several canonical models exist -- for example, the Preece-Baines and Jolicoeur models \cite{JPPS} -- they are specialized to human development and as such have an outsized number of parameters. Without knowing in detail the particular characteristics of dragon growth over a wider population, we will not speculate about the phases of dragon growth as might motivate a more complex model. Instead, we use the naive model of a sigmoid function. As such, the form of $L(t)$ becomes (for some as yet undetermined $a$ and $b$):

\[\dfrac{A}{1+ae^{-bt}}\]

We then identify the width as:
\[W(t) = \dfrac{1}{L(t)}\int_{0}^{L(t)}W(l)dl\]
Practically, this indicates the discrete summation we implemented over the cross-sectional areas of various lizard species and various images of Drogon, making $W(t)$ the average width at a given time $t$ -- so while $W(t)$ is not exact, the principle still stands.

We motivate this explicit definition of $W(t)$ as an attempt to express the mass, $M(t)$. Using the assumptions laid out in Section 5.1, the argument for an approximate definition of the mass in terms of $W(t)$ and $L(t)$ as known quantities becomes relatively straightforward.

By assumption one, we have,

\[W(t) \;\si{m} \propto H(t) \;\si{m} \] 
\[W(t)^2 \;\si{m^2} \propto W(t) \cdot H(t) \;\si{m^2}\]

By assumption two, 

\[\dfrac{1}{L(t) \;\si{m}}\int_{0 \;\si{m}}^{L(t) \;\si{m}}C(l) \;\si{m^2} dl \propto W(t) \cdot H(t) \;\si{m^2}\]

\[\dfrac{1}{L(t)}\int_{0}^{L(t)}C(l)dl \;\si{m^2} \propto W(t)^2 \;\si{m^2}\]

\[\int_{0}^{L(t)}C(l)dl \;\si{m^3} \propto L(t) \cdot W(t)^2 \;\si{m^3}\]

\[\rho \;\si{kg/m^3} \int_{0}^{L(t)}C(l)dl \;\si{m^3} \propto \rho \;\si{kg/m^3} \cdot L(t) \cdot W(t)^2 \;\si{m^3}\]

By assumption four,

\[\rho\int_{0}^{L(t)}C(l)dl \;\si{kg} \approx M(t) \;\si{kg}\]

\[M(t) \;\si{kg} \appropto \rho \cdot L(t) \cdot W(t)^2 \;\si{kg}\]

This gives us a relatively explicit statement of $M(t)$.

\[M(t) \;\si{kg} \appropto \rho \;\si{kg/m^3} \;L(t) \;\si{m} \left(\dfrac{1}{L(t) \;\si{m}} \int_{0 \;\si{m}}^{L(t) \;\si{m}}W(l) \;\si{m}\;dl \right)^2\]

\[M(t) \;\si{kg} \appropto \dfrac{\rho}{L(t)} \left(\int_{0}^{L(t)}W(l)dl \right)^2 \;\si{kg}\]

\subsection{Modeling}
\hspace{1em}
We will also add as an addendum that

\[V(t) \;\si{m^3} \approx V_T(t) \;\si{m^3} = \int_{0 \;\si{m}}^{L(t \;\si{m})}C(l) \;\si{m^2}dl\]

such that the volume of the torso as the cross-sectional area summed over the length of the dragon is approximately equal to the exact volume of the dragon, which makes implicit an assumption that $\rho V_T(t) \appropto M(t)$. More appropriately, it alters our definition of $\rho$ such that $\rho$ is a constant that when multiplied by $V_T(t)$ is the mass of $M(t)$.

We initially began the simulation of the volume based on our own measurements of the length as obtained in Section 4, Data Collection. We then found a 3D model of a dragon online, using the software Blender to simulate the dragon's volume and the length of its head. We then scaled the measurements based on our directly observed data over time. Finally, we obtained the mass by multiplying by $\rho$, which we took as its upper bound, the density of the water, per assumption 3.

Letting $L_{H}^{M}$ signify the constant length of the model's snout and $V_{T}^{M}$ the constant volume of the whole model, the equation for mass follows:

\[M(t) = \dfrac{\rho V_{T}^{M}}{\left(\dfrac{L_{H}^{M}}{L_H(t)}\right)^3}\]
\[M(t) = \rho V_{T}^{M}\left(\dfrac{L_H(t)}{L_{H}^{M}}\right)^3\]

This approach was ultimately abandoned due to two implied assumptions. First, to apply a scaling as we did in the denominator of the original equation, we must assume $L_H(t)^3 \propto V(t)$. By the final equivalence in 5.2,

\[V(t) \approx V_T(t) = V_{T}^{M}\left(\dfrac{L_H(t)}{L_{H}^{M}}\right)^3\]

Second, we must assume that various dimensions of the dragon grow proportionally to one another -- in particular, that $L_H(t) \propto L(t)$, which is deeply problematic. Independent of the many issues with the approach, though, it still remains a useful heuristic for obtaining a workable mass, acting as a satisfactory upper bound, if a crude one. As such, we include the mass data from this na\"{i}ve approach below.

\begin{table}[h]
\begin{tabular}{|l|l|l|l|l|l|l|l|}
\hline
Age (years) & 0    & 0.5  & 2     & 3       & 4       & 5        & 6+        \\ \hline
Mass (kg)   & 2.60 & 3.70 & 98.38 & 1947.54 & 7544.14 & 90089.83 & 251328.27 \\ \hline
\end{tabular}
\end{table}

\subsection{Extrapolated Energetic Consumption}
\hspace{1em}
After obtaining the mass of the dragons, we attempted to use metabolic rates to provide a satisfactory conclusion to our central question concerning the feasibility of a dragon's existence in our reality.

We use the metabolic equation of Raveling and LeFebvre's \emph{Energy Metabolism and the Theoretical Flight Range of Birds} \cite{metabolism} to obtain the basal metabolic rate of a dragon in \si{kcal/hr}. Rewriting the equation they provide for \si{kcal/hr} with respect to mass in kilograms (originally sourced from a 1961 article by King and Farner) in terms of our notation,

\[\log(R_B(t)) = \log(3.1) + 0.744\log(M(t)) + 0.074\]
\[R_B(t) = \left(3.1 \cdot 10^{0.074}\right)\left(M(t)^{0.744}\right)\]

According to the same paper by Raveling and LeFebvre, in general, the energetic expenditure during migratory flight for most birds (extrapolating migratory flight as a behaviour to characteristic draconic flight, e.g. the flight of Balerion to Valyria in Fire and Blood \cite{fireandblood}) is about twelve times the basal metabolic rate, such that $12R_B(t) = R_F(t)$ (implying that $\frac{3}{2} R_B(t) = R_S(t)$). This gives us an equation for total energetic consumption $C(t)$ in \si{kcal/day} of flight as follows:

\[C(t) = 24\left(P_{F}R_{F}\left(t\right)+\left(1-P_{F}\right)R_{S}\left(t\right)\right)\]

Or, alternatively:

\[C(t) = 24\left(12P_{F}R_{B}\left(t\right)+\dfrac{3}{2}\left(1-P_{F}\right)R_{B}\left(t\right)\right)\]

Using the projected mass in Section 5.3.1 from ages 6 and 7 of 251328.27 kg, this expression for energetic consumption gives us in the worst-case scenario of $P_F:=1$ a daily consumption in kilocalories of 11,028,939. Using a far more realistic $P_F$ value of $\frac{1}{3}$, we arrive at a projected consumption of 4,595,391 \si{kcal/day} around the ages of 6 and 7, which squares relatively well with our understanding of a dragon's daily food intake.

 \subsection{FM Summary}
\hspace{1em}
Overall, we simulated the volume of the dragon by using scaled measurements to obtain a general number across stages and then utilized these volumes and the metabolic rates to get the energy consumption of dragons at the ages of around 6 or 7.

From the numbers, we can safely say that the forward model proves that the existence of a dragon is possible simply from a energy intake view as an average, 63.5 sheep would only provide 21 kg of edible meat which is 54180 kcals, and the required daily intake of over 4 million kcals would require the dragon to be able to consume around 85 sheep per day, which is entirely possible.

\section{Backward Model}

\subsection{Introduction}
\hspace{1em}
Differing from the forward model introduced in the previous section, we build another model to analyze dragon's energy consumption and expenditure in a backward way. Given the literal data of food a dragon consumed, we calculate the amount of energy the diet received by the dragon, and use this result to explore the plausibility of dragon existence through a empirical examination on energy spent on flying, fire-breathing, etc.

Based on \textit{A Song of Ice and Fire} by George R.R. Martin, sheep are a Dragon's favorite food, so the energy consumption calculation will be simplified to be derived only from sheep.

\subsection{Assumptions}
\hspace{1em}
During the process of modeling, there is a number of assumptions made as to the following:
\begin{enumerate}
    \item Magic:
    \begin{enumerate}
        \item The initial stages of modeling is built without magic
	BUT later the interference of magic may be key to explain things like “conservation of energy”.
	\end{enumerate}
	\item Wisdom:
	\begin{enumerate}
	    \item the wisdom of dragons is set to be lower than humans, which means dragons do not have their own civilizations or abilities such as producing their own tools.
    \end{enumerate}
\end{enumerate}

\subsection{Data for Energy Consumption}
	
First and foremost, we attempt to answer the first question: “What is the data of dragons from fantasy literature?”

The data we consider includes mass, age, diet, energy consumption, reproduction, habits, etc.

\begin{enumerate}
\item Energy Conversion Efficiency
As sheep are dragons’ favorite source of food, according to George R.R. Martin in A Song of Ice and Fire, we will compute the energy conversion efficiency of dragons focusing on sheep.

\item Sheep Raising 

100 acres of land can support about 280 female sheep \cite{sheep_raising}.


\item Market Lambs are usually 6-8 months old
\item Weighs 140 lbs (63.50 kg), yielding 46 to 49 pounds (20.87 to 22.22kg) of edible lean retail lamb cuts, semi-boneless.

\item Around 15 kilograms of feed turns to 1 kilogram of goat/sheep meat. Energy feed conversion efficiency of goat/sheep meat = 4.4\% \cite{diet}. 

\item Energy Cost of Flying

For calculation of energy consumption, the formula that will be using is based on previous section which is \cite{metabolism}

\begin{equation}
    C = 24\left(12P_{F}R_{B}\left(t\right)+\dfrac{3}{2}\left(1-P_{F}\right)R_{B}\left(t\right)\right)
\end{equation}

where the $R_B$ is given as 
\begin{align}
    R_B = \left(3.1 \cdot 10^{0.074}\right)\left(M^{0.744}\right)
\end{align}

\item Diet 

 Komodo dragon's diet behavior will be used as base model for the dragon diet behavior. Based on a Natgeo article as well data from San Diego zoo, Komodo dragon can eat up to 80\% of its weight and has digestive pause ranging from 3-6 days. Thus, we can assume that Komodo dragon can eat prey about 13\% of its weight daily \cite{natgeokomodo}\cite{zookomodo}.  
\end{enumerate}

\subsection{Modeling}

After collecting data from both fantasies and our mother nature, finally, we come to the exciting stage of modelings. 

What are the ecological impacts and requirements of the dragons? 

To support such mighty creatures’ living, it takes serious ecological impacts and requirements. As suggested earlier, consider sheep to be the only food source for dragons, a regular dragon eats a lamb a day. Whereas a 2000-ton dragon eats 15 sheep a day.

As a regular dragon weighed about 490 kg if we based off Komodo dragon's daily diet assuming the sheep's weight is at 63.5 kg. 

From this, we get to do a simple construction of the formula describing a dragon’s weight to the number of sheep it requires to consume:

\begin{equation}
    M_r = \dfrac{63.5}{0.13}N
\end{equation}

Now, a further question is raised: “How many sheep support a sheep to be consumed in the ecosystem?” If a 100 acre land can raise 280 female sheep, thus the quick model describing the relationship between the number of sheep support and number of sheep that need to support the ecosystem is 

\begin{equation}
    N_{sheep support}= \frac{N}{280}
\end{equation}
 
In addition, the required elements for the land, mainly the size, to support these lambs to grow naturally / by humans? For which we also need to decide a place to meet all the traits (grass, lambs … ) e.g. Yellowstone 

Now, taking consideration of sheep per dragon, the sample size of sheep to sustain sheep eaten by dragons (first ignoring the natural growth of dragons), test whether it is possible to find a place for such ecosystem to exist.

Consumption vs. Mass on average per day
This part of modeling attempts to establish a relationship between dragons (size, mass, age) and their energy expenditure.



\begin{equation}
m = (Variable) \cdot t
\end{equation}

As mentioned in previous part, dragons’ energy consumption/expenditure is made up of two main sources, flying and fire-breathing. Thus, the total energy expenditure will be considered of these two plus “others” which include any other energy outputs. The determining factor is m as of mass in kg. 

\begin{equation}
E_{Expenditure}(m) = E_{Flying}(m) + E_{Fire-Breathing}(m) + E_{Others}(m)
\end{equation}

The energy required to cook the lamb based on 3000 W oven for 120 minutes will be about 21600000 J. This will be the approximation for the energy consumption for fire breathing for a dragon per sheep. 

$E_{comsumption}(E_{expenditure}(m)) = k E_{expenditure}(m)$

"The dragons this year are the size of 747s," Matt Shakman, one of next season's four directors, told Entertainment Weekly. (That's about 230 feet long and 210 feet in wingspan, according to the mag.) "Drogon is the biggest of the bunch—his flame is 30-feet in diameter!" Shakman added.

Nettles tamed Sheepstealer by bringing him a freshly slaughtered lamb every day until he grew accustomed to her presence -- in-universe manuscript art.

The energy needed to fly daily for a dragon is 

\begin{align}
E_{Flying} &=4184\cdot C \\ &= 4184 \cdot 24 \left(12P_{F}R_{B}\left(t\right)+\dfrac{3}{2}\left(1-P_{F}\right)R_{B}\left(t\right)\right)\\ 
\end{align}

where the $E_{Flying}$in J.  

Let n be the number of sheep a dragon eats per day
\begin{align}
E_{Fire} &= 21600000 \cdot n \\
&= 21600000\cdot20 \\
&= 4.32\cdot10^8 J    
\end{align}

If a dragon has $m = 2,000,000 kg = 2\cdot10^6 kg$ and assume it can travel as long as 8 hours or $P_F = \frac{1}{3}$, so the energy for flying will be about $8.99\cdot10^{10}$ J.

From Lambs, we can get 294 calories per 100 g, which is approximately 3000 calories / kg = 12552 J/kg = 941400 J per lamb.

Thus, the number of lamb needed to sustain 2000-ton dragon daily is about 95571, which translates to about 34200 acres of land or about 70\% Seattle area. Considering the size of the dragon, it makes sense that it will need such large amount of food to sustain its daily activities.  

Thus, we only need to discuss the enormous energy consumption that is needed by the dragon. The total energy provided by 20 lambs is about $1.8\cdot10^7$, which means there are a more efficient energy conversion inside the body of the dragon that is breaking the known thermodynamic laws. We can further explain with the help of "scientific" magic. We can think that there is a presence of dragon crystal that is able to convert the energy from sheep meat into stronger energy type. As a result, dragon can not exist in the real world due to presence inexplicable dragon crystal.


 \subsection{BM Summary}


The backward model disproves the existence of a dragon based on the fact that it requires energy crystal in order to convert the energy from food into fire breathing ability, which consumes more energy than energy that dragon gets from eating sheep.

\section{Conclusion}
\hspace{1em}
This paper introduces two strategies to analyze the possible existence of Game of Thrones dragons. The forward model supports the existence of dragons while the backward model disproves the hypothesis. We ultimately conclude that dragons cannot exist as described in the relevant literature. The signature behavior of the fantasy dragon, fire breath, is ultimately too costly in energy to be plausibly supported by a reasonable environment. However, having only accounted for flight via metabolic rate calculations in the forward model, we conclude that a dragon without fire breath could be supported.

In conclusion, if dragons had naturally evolved as flying animals unable to breathe fire and cause mass destruction, they would theoretically be able to live and fly in the real world as illustrated by the forward model. However, due to the fantastical characteristics of dragons, they would be unable to exist on Earth. Their signature ability, fire breath, simply costs far more energy than it is reasonable for any animal to expend and would require impossibly vast swathes of land.

\section{Acknowledgements}
\hspace{1em}
The inspiration for this research came from Zheng Cao, Wenyu Du, and Zhuorui He, all undergraduates at University of Washington, in 2021. Wenyu and Zhuorui helped construct the draft outline as well as offering some insightful suggestions to the paper. while Zheng continued leading the research with Brody, Jiayi, Mark, and Vinsensius.

\end{document}